\documentclass{article}

\usepackage{arxiv}

\usepackage[utf8]{inputenc} 
\usepackage[T1]{fontenc}    
\usepackage{hyperref}       
\usepackage{url}            
\usepackage{booktabs}       
\usepackage{amsfonts}       
\usepackage{nicefrac}       
\usepackage{microtype}      
\usepackage{lipsum}

\usepackage{epsfig}

\newtheorem{theorem}{Theorem}

\newtheorem{remark}{Remark}

\usepackage{verbatim}
\title{Non-integrability of  the axisymmetric Bianchi IX cosmological model via Differential Galois Theory}

\author{Primitivo ACOSTA-HUM\'ANEZ \\
Instituto de Matem\'atica \& Escuela de Matem\'atica\\
\href{www.uasd.edu.do}{Universidad Autonoma de Santo Domingo},\\
Dominican Republic \\
  \texttt{pacosta-humanez@uasd.edu.do} \\
   \And
Juan J.~MORALES-RUIZ \\
~Depto. de Matem\'atica Aplicada\\ \href{www.upm.es}{Universidad Polit\'ecnica de Madrid}\\ Madrid, Spain\\
    \texttt{juan.morales-ruiz@upm.es} \\
   \AND
  Teresinha J. STUCHI \\
  Instituto de F\'isica\\ \href{www.ufrj.br}{Universidade Federale de Rio de Janeiro}\\ Rio de Janeiro, Brazil\\
   \texttt{tstuchi@if.ufjr.br} \\
}

\begin{document}
\maketitle

\begin{abstract}
We investigate the integrability of an anisotropic universe with matter and cosmological constant formulated as Bianchi IX models. The presence of the cosmological constant causes the existence of a critical point in the finite part of the phase space. The separatrix associated to this Einstein's static universe is entirely contained in an invariant isotropic plane forgetting the singularity at the origin. This invariant plane of isotropy is an integrable sub-space of the Taub type. In this paper we analyse the differential Galois group of the second order variational equations to this plane in order to apply the integrability theorem of the second author with Ramis and Sim\' o.  The main result is that the model is non-integrable by meromorphic functions.
\end{abstract}

\section*{Introduction}

Belinskii, Khalatnikov and Lifishitz \cite{belinskii2} started the question of chaotic behaviour of general Bianchi IX models in Relativistic Cosmology. The interest in the chaoticity (or not) of Bianchi IX models has been mainly focused on the Mixmaster case. Vacuum Bianchi IX models with three scale factors was taken by Misner \cite{misner}).   Also  \cite{bogo,bono} studied the three dimensional
Bianchi IX model from the real dynamics point of view. The absence of an invariant (or topological) characterisation of chaos in the model,  i.e., standard chaotic indicators such as Lyapunov exponents being coordinate dependent and therefore questionable \cite{matsas} ,\cite{calzeta}. For discussions of the issue of chaotic dynamics on these models we refer to the works of Barrow \cite{barrow}. Cornish and Levin   \cite{cornish}  proposed to quantify chaos in the Mixmaster universe by calculating the dimensions of fractal basin boundaries in initial-conditions sets for the full dynamics. Calzeta \cite{calzeta} and El Hasi proved integrability under adiabatic approximation, but the tori breaked up without the adiabatic terms.

We recall that the conjunction of a cosmological constant and anisotropy implies in the existence of a critical point, the Einstein universe, in the finite region of phase space. This point is linearly a saddle center that implies the existence of a centre manifold tangent to the linear periodic orbits around the associated energy $E$. From this centre manifold the stable and unstable manifolds which are topologically cylinders  carry the dynamics away and to the inflationary region.  This was the object of the work of H.P. Oliveira et al. (see \cite{oliveira}) for the axisymmetric Bianchi IX model. This is the model we study in this paper. In \cite{oliveira} this problem was numerically shown non integrable; these cylinders, when the full dynamics is considered,  intercept each other in homoclinic and heteroclinic connections causing a chaotic dynamics. This has as consequence a chaotic escape in the inflationary regime. In \cite{mora2}, see also \cite{mo}, is showed  that the three degrees of freedom  Bianchi IX models, without cosmological constant is no-integrable by rational functions. The axisymmetric version has so far no analytic proof of non integrability.

In this paper we show that the two degrees of freedom axisymmetric Bianchi IX model with cosmological constant is non integrable by meromorphic functions for generic values of the parameters and by assuming a natural conjecture. To obtain this result we apply a joint theorem with Ramis and Sim\'o in \cite{morasi}, where it is given a necessary condition for integrability of complex analytical Hamiltonian systems by means of meromorphic first integrals, the identity components of the Galois groups of the higher order variational equations must be commutative, see Theorem \ref{tmorasi} for a precise statement. This result was an extension to higher order variational equations of a previous joint theorem with Ramis and the second author in \cite{mora1}, see Theorem \ref{tmora1} in Appendix \ref{sec:DGT}.\medskip

 In section 1 we present  the axisymmetric model.  In  section 2  the  calculation of Galois groups for first and second  variational equations are made. In section 3 we illustrate
the behaviour of the system with Poincar\'e sections.  An overview about Differential Galois Theory and non-integrability of Hamiltonian systems are presented in \ref{sec:DGT}. Finally, application of Kovacic algorithm to compute the solutions and Galois of the first variational equation is presented in \ref{sec:kovacic}.

\section{The axisymmetric  Bianchi IX Inflationary model} 

The axisymmetric Bianchi IX models, see \cite{stuchi},  derived from the full three scales model, $H(A,B,C,P_A,P_B,P_C)$   considering  $B=C$, $\dot B=\dot C$,
 with cosmological constant  $\lambda$ is:

\begin{equation}\label{eq:log12}
{H}(A,B,P_A,P_B)=\frac{P_AP_B}{4B}-\frac{AP_A^2}{8B^2}+2A-\frac{A^3}{2B^2}-2\lambda AB^2
\end{equation}

 It obeys the Einstein equation :
  $$G^{\mu\nu}-\lambda g^{\mu\nu}= T^{\mu\nu}.$$ and
 the line element is obviously $ds^2=dt^2 -A({t})^2(\omega^1) ^2-B(t)^2((\omega^2)^2+(\omega^3)^2)$.


The introduction of the cosmological constant is one of the cornerstones in the paradigm of inflation.   The inflationary constant  plays a fundamental role in the gravitational dynamics model, by
inducing an expansion of the scales towards the De Sitter configuration.\\

The equations of motion corresponding to Eq. \eqref{eq:log12} are given by the Hamiltonian system

\begin{equation}\label{eq:two}
\begin{aligned}
\dot A=&\frac{P_B}{4B}- \frac{AP_A}{4 B^2}\\ \dot B=&\frac{P_A}{4B}\\
\dot{P_A}=&\frac{P_A^2}{8B^2 }-2+\frac{3A^2}{2B^2}+2\lambda B^2\\
\dot P_B=&\frac{P_AP_B}{4B^2}-\frac{AP_A^2}{4B^3}-\frac{A^3}{B^3} +4\lambda AB
\end{aligned} \end{equation}

This system,  Eq. \eqref{eq:two}, has a critical point $\mathcal E$ in the finite region of the phase space with
coordinates
$$\mathcal E: A_0=B_0=\frac{1}{\sqrt{4\lambda}}, \quad P_A=P_B=0, $$

\noindent with associated energy $E_0=1/\sqrt{4\lambda}$. This point corresponds to
the static Einstein universe. From the point $\mathcal E$ emanate the stable and
unstable manifolds which goes to infinity  to the right of $\mathcal E$ and  in the left
the unstable and stable branch meet at the origin. The origin is a singular
point: The Friedmann universe

Furthermore system \eqref{eq:two} has an invariant plane $\Gamma$: $A=B$, $P_B=2P_A$
that gives us a family of integral curves, the Taub family, see \cite{taub}, defined by the one degree of freedom Hamiltonian

\begin{equation}\label{cubham}
E=\frac{p_x^2}{8x}+\frac 32 x-2\lambda x^3,
\end{equation}
where $x:=A=B$ and $p_x:=P_A=\frac{1}{2}P_B=4x\dot x$.

\section{Computation of the variational equations}

As we want to compute the variational equations along the separatrix of the Taub family as
a second order variational equation, it is convenient to write de hamiltonian equations as
second order differential equations
\begin{equation}\label{hebi}
    \begin{array}{ccl}
 \ddot {A}&=&\frac{A}{2B^2} -\frac{5A^3}{8B^4} +\frac{A \dot B^2}{2B^2}+\frac{\lambda A}{2} -\frac{\dot A \dot B}{B}\\
\ddot B&=&\frac{3A^2}{8B^3}-\frac{1}{2B}  +\frac{\lambda B}{   2 } -\frac{\dot B^2}{2B}
\end{array}
\end{equation}
The invariant plane $\Gamma$ is given in these coordinates by $A=B$, $\dot A= \dot B$, and the Taub family along this plane is parametrized by coordinates $x:=A=B$ and  $\dot x=\dot A=\dot B$
and given by
\begin{equation}\label{eqinvplane}\ddot x=-\frac{1}{8x} +\frac{\lambda x}{2}  -\frac{(\dot x)^2}{2 x}\end{equation}

The first integral  of the Hamiltonian restricted to the invariant plane $\Gamma$ gives the integral
curves of the above equation which in the coordinates $(x,\dot x)$, through Eq. \eqref{cubham} because $y=4x\dot x$,
 are expressed as the family of cubics
\begin{equation}\label{eqalpha}
\Gamma_E:6 x \dot x^2+\frac{3x}{2}-2\lambda x^3-E=0,
\end{equation}
for generic values of $\lambda$ and $E$. \subsection{Variational Equation}
 The variational equations of the Hamiltonian system \eqref{hebi} along a particular solution of $\Gamma$ is given by  the equations
\begin{equation}\label{bnveg}
\begin{array}{ccl}
\ddot\psi_A&=& \left(\frac{1}{2B^2} -\frac{15A^2}{8 B^4} +\frac{\dot B^2}{2B^2}+\frac{\lambda}{2}\right)\psi_A-\frac{\dot B}{B}\dot \psi_B
 \\&& +\left(\frac{2\lambda A}{B}-\frac{3\dot A\dot B}{B^2}-\frac{4\ddot A}{B}+\frac{A\dot B^2}{B^3}+\frac{A}{B^3}\right)\psi_B+
 \left (\frac{A \dot B}{B^2}-\frac{\dot A}{B}\right)\dot\psi_B\\
\ddot \psi_B&=&\frac{3A}{4B^3}\psi_A+(\frac{-9A^2}{8B^4}+\frac{1}{2B^2}+\frac{\lambda}{2}+\frac{\dot B^2}{2B^2})\psi_B-\frac{\dot B^2}{B^2}\dot \psi_B
\end{array}
\end{equation}

\noindent where the coefficients are restricted to the integral curve $\Gamma$: $A=A(t)$, $\dot A=\dot A(t)$,
$ B=B(t)$ and $\dot B=\dot B(t)$.

For $\Gamma=\Gamma_E$, being $\psi=\psi_A$, $\zeta=\psi_B$, $x=A=B$, $\dot x=\dot A=\dot B$, we obtain the VE (variational equations) along the Taub solution,
\begin{equation}\label{FVEG}
\begin{array}{lll}
\ddot \psi&=&-\frac{\dot x}{x}\dot\psi+\left(\frac{1}{2}\left(\frac{\dot x}{x}\right)^2-\frac{11}{8x^2}+\frac{\lambda}{2}\right)\psi\\&& -\left(\frac{4\ddot x}{x}+2\left(\frac{\dot x}{x}\right)^2-\frac{1}{x^2}-2\lambda\right)\zeta\\
\ddot \zeta&=&-\frac{\dot x}{x}\dot\zeta-\left(\frac{4\ddot x}{x}+\frac{3}{2}\left(\frac{\dot x}{x}\right)^2+\frac{9}{8x^2}-\frac{5}{2}\lambda\right)\zeta +\frac{3}{4x^2}\psi
\end{array}
\end{equation}
We set the following notations that will be useful along the rest of the paper:
\begin{equation}\label{cubpols}
  C_1(x) := 4\lambda x^3-3x+2E,\quad C_2(x):=35\lambda x^3+210x+4E
\end{equation}
Now in order to obtain a differential equation over the Riemann sphere ${\bf P}^1$
(i.e. with rational coefficients), it is convenient to apply the so-called \emph{Hamiltonian algebrization procedure}, see \cite{acbl,amw}. This procedure starts with a \emph{rational Hamiltonian change of variable}, that is,  making a change of variables $t\mapsto x$ with $x:=x(t)$ such that $(\dot x)^2=\alpha(x)\in\mathbb{C}(x)$, for instance, $\dot x=\pm\sqrt{\alpha(x)}$. This change of variable is obtained by the using of equation of $\Gamma_E$, equation \eqref{eqalpha}, and then we have
$$\alpha=\frac{C_1(x)}{12x},\quad \psi'=\sqrt{\alpha}\widehat{\psi}',\quad \zeta'=\sqrt{\alpha}\widehat{\zeta}',$$ where $C_1(x)$ is given in Eq. \eqref{cubpols}. Applying Hamiltonian algebrization procedure, where Eq. \eqref{eqinvplane} provides $\ddot x$, we obtain 
the algebraic variational equations (AVE)
\begin{equation}\label{eq:ave}
    \begin{array}{ccc}
  \widehat{\psi}''=-\frac{8\lambda x^3-3x+E}{xC_1(x)}\widehat{\psi}'+\frac{8\lambda x^3-18x+E}{x^2C_1(x)}\widehat{\psi}+\frac{18}{xC_1(x)}\widehat{\zeta}\\ \\
   \widehat{\zeta}''=-\frac{8\lambda x^3-3x+E}{xC_1(x)}\widehat{\zeta}'+\frac{8\lambda x^3-9x+E}{x^2C_1(x)}\widehat{\zeta}+\frac{9}{xC_1(x)}\widehat{\psi}
     \end{array},
\end{equation} 
for generic values of $\lambda$ and $E$, which correspond to the algebraic form of variational equations given in equation \eqref{FVEG}.

After a hard work combining computations by hand with Maple, we solve the algebraic form of the variational equation \eqref{eq:ave}. Therefore one basis of solutions of the variational equation \eqref{eq:ave} is given by 
$$\mathcal{B}_{fve}=\left\{\begin{pmatrix}\widehat{\psi}_1(x)\\ \widehat{\zeta}_1(x)\end{pmatrix},\begin{pmatrix}\widehat{\psi}_2(x)\\ \widehat{\zeta}_2(x)\end{pmatrix},\begin{pmatrix}\widehat{\psi}_3(x)\\ \widehat{\zeta}_3(x)\end{pmatrix},\begin{pmatrix}\widehat{\psi}_4(x)\\ \widehat{\zeta}_4(x)\end{pmatrix}\right\},$$ where $\widehat{\psi_i}(x),\widehat{\zeta_i}(x)$ are given by
\begin{equation}\label{eq:foave}
\begin{array}{l}\widehat{\psi}_1(x)=\widehat{\zeta}_1(x)=\sqrt{\frac{C_1(x)}{x}}=\sqrt{12\alpha}\\ 
\widehat{\psi}_2(x)=\widehat{\zeta}_2(x)=\sqrt{\frac{C_1(x)}{x}}\int\sqrt{\frac{{x}}{C_1(x)^3}}dx=\widehat{\psi}_1(x)\int\frac{dx}{x\widehat{\psi}_1(x)^3}\\ \widehat{\psi}_3(x)=-2\widehat{\zeta}_3(x)=\frac{2}{3}\sqrt{\frac{C_2(x)}{x}}\exp\left(-i\int\frac{ \sqrt{486(E^2\lambda+2450)}dx}{C_2(x)\widehat{\psi}_1(x)}\right)\\
\widehat{\psi}_4(x)=-2\widehat{\zeta}_4(x)=\frac{2}{3}\sqrt{\frac{C_2(x)}{x}}\exp\left(i\int\frac{ \sqrt{486(E^2\lambda+2450)}dx}{C_2(x)\widehat{\psi}_1(x)}\right)\end{array}\end{equation}
 Thus, the basis field for Eq. \eqref{eq:ave} is $K=\mathbb{C}(x)$ and the Picard-Vessiot extension is $L=K(\widehat{\psi}_i,\widehat{\psi}_i')$, $i=1,2,3,4$. For instance, the differential Galois group of Eq. \eqref{eq:ave} is a subgroup of $GL(4,\mathbb{C})$.

\subsection{Normal Variational Equation}
We observe that in  $\Gamma_E$,  the equation (\ref{FVEG})  is satisfied by the  change of variable   $\chi:=\psi-\zeta$, leading to the normal variational equation:

\begin{equation}\label{nveg}\ddot \chi+\left(\frac{\dot x}{x}\right)\dot \chi+\left(\frac{17}{8x^2}-\frac{\dot x^2}{2x^2}-\frac{1}{2}\lambda\right)\chi=0
\end{equation}
being $x=x(t)$, $\dot x=\dot x(t) $ the Taub solution in  $\Gamma_E$. Recalling that $$\alpha=\frac{C_1(x)}{12x},$$ we obtain the algebraic form of VE corresponding to equation \eqref{nveg}, through the change of variables $(x(t),\dot x(t),\chi,\dot \chi, \ddot \chi)\mapsto (x,\sqrt{\alpha},\eta,\sqrt{\alpha}\eta,\frac{1}{2}\alpha'\eta'+\alpha\eta'')$, which is the AVE given by 

\begin{equation}\label{eq2o}
\begin{array}{l}    \eta''+  p(x)\eta'+q(x)\eta=0,\\  \\p(x)=\frac{8\lambda x^3-3x+E}{xC_1(x)},\, q(x)=\frac{-8\lambda x^3+27x-E}{x^2C_1(x)},
\end{array}
\end{equation}
being $C_1(x)$ provided in Eq. \eqref{cubpols}. Now, the invariant normal form of the equation \eqref{eq2o} can be obtained as
$$
\frac{d^2\xi}{dx^2}-g(x)\xi=0,\, g(x)=\frac{p^2(x)}{4} +\frac{1}{2}\frac{dp(x)}{dx} - q,\, \xi=\eta\cdot\exp\left(\frac{1}{2}\int_{x_0}^xp(t)dt\right)
$$

Thus, the invariant normal form of equation \eqref{eq2o} corresponds to

\begin{equation}\label{eq:ggen}
\begin{array}{l}
\frac{d^2\xi}{dx^2}=g(x)\xi,\,\, \xi=\eta\sqrt[4]{xC_1(x)}\\ \\
  g(x)=\frac{128\lambda^2 x^6-552\lambda x^4+128E\lambda x^3+315x^2-222Ex+5E^2}{4x^2C_1(x)^2}.
  \end{array}
\end{equation}
Because the Galois group of the normal variational equation is a subgroup of $SL(2,\mathbb{C})$, we can apply the known results about the classification of subgroups of $SL(2,\mathbb{C})$, see \cite{kovacic:1986} and references therein.\\

We remark that avoiding the trivial cases $\lambda=0$ ($C_1(x)=0$ reduces to linear equation) and $E=0$ ($C_1(x)=0$ reduces to a quadratic equation), the discriminant for the cubic equation $C_1(x)=0$ is $$\Delta=-\frac{27(4E^2\lambda-1)}{16\lambda^3},\quad C_1(x)=4\lambda x^3-3x+2E=0$$ Moreover, for $\lambda=(4E^2)^{-1}$ the discriminant vanishes ($\Lambda=0)$.  The discriminant $\Lambda$ is an important element for the galoisian analysis of the above differential equation, that is, equation \eqref{eq:ggen}. Furthermore, the differential equations \eqref{bnveg}, \eqref{eq:ave}, \eqref{nveg}, \eqref{eq2o} and  \eqref{eq:ggen} over the invariant plane $\Gamma_E$ have the same connected identity component of the Galois group because they are linked through change of variables (algebraic) that allows to preserve the connected identity component of the Galois group, see \cite{mo}. 

We recall that a Fuchsian differential equation is a linear homogeneous ordinary differential equation with analytic coefficients in the complex domain whose singular points are all regular singular points, see \cite{abst}. 
\begin{remark}
We know that if the first variational equation is fuchsian, then the next variational equations are fuchsian too. Thus, the solutions given by Eq. \eqref{eq:foave} have no exponential behaviour around the singular points. In fact, it is possible to verify the existence of enough logarithms in the integrand in the exponentials of $\widehat{\psi_i}$, $i=3,4$, in Eq. \eqref{eq:foave}.
\end{remark}

Therefore, the differential equations \eqref{eq:ave}, \eqref{nveg}, \eqref{eq2o} and  \eqref{eq:ggen} over the invariant plane have the same set of singularities which is, using the Cardano's formula for the cubic equation, given by the set of poles of $g$ including $\infty$ as follows: $\Gamma=\{0,\infty,\rho_0,\rho_1,\rho_2\}$, where $$ \rho_j=\nu^j\sqrt[3]{-\frac{E}{4\lambda}+\frac{i\sqrt{3\Delta}}{18}}+\nu^{2j}\sqrt[3]{-\frac{E}{4\lambda}-\frac{i\sqrt{3\Delta}}{18}},\, \nu^3=1,\, j=0,1,2.$$ For suitability, we write the cube roots as follows:
    $$\rho_j=\nu^j\frac{\kappa}{2\lambda}+\nu^{2j}\frac{1}{2\kappa},\, \kappa=\sqrt[3]{\left(-2E+\sqrt{\frac{4E^2\lambda-1}{\lambda}}\right)\lambda^2},\, \nu^3=1,\, j=0,1,2.$$
These five singularities are of the regular type for $\Lambda\neq 0$, which implies that these differential equations are fuchsian differential equations with five regular singularities. If either $E=0$, $\lambda=0$ or $\Lambda=0$, at least one singularity is of irregular type.

To obtain the solutions of Eq. \eqref{eq2o} and Eq. \eqref{eq:ggen}, we consider $\Delta\neq 0$, $E\neq 0$ and $\lambda\neq 0$. Due to $\chi=\psi-\zeta$ we have $\eta=\widehat{\psi}-\widehat{\zeta}$. Thus, setting $$\mu=i\sqrt{486(E^2\lambda+2450)}\neq 0$$ we have the following solutions for $\eta$ in Eq. \eqref{eq2o}, according to Eq.  \eqref{eq:foave}:
\begin{equation}\label{eq:sed2o}
    \begin{array}{l}
 \eta_1(x)=\widehat{\psi}_4(x)-\widehat{\zeta}_4(x)=\sqrt{\frac{C_2(x)}{x}}\exp\left(\mu\int\frac{\sqrt{x}dx}{C_2(x)\sqrt{C_1(x)}}\right)\\
 \eta_2(x)=\widehat{\psi}_3(x)-\widehat{\zeta}_3(x)=\sqrt{\frac{C_2(x)}{x}}\exp\left(-\mu\int\frac{\sqrt{x}dx}{C_2(x)\sqrt{C_1(x)}}\right)\\
 \eta_3(x)=\widehat{\psi}_2(x)-\widehat{\zeta}_2(x)=0\\
  \eta_4(x)=\widehat{\psi}_1(x)-\widehat{\zeta}_1(x)=0
\end{array}
\end{equation}
For suitability we reordered the basis of solutions of Eq. \eqref{eq:foave} to obtain as basis of solutions of Eq. \eqref{eq2o} the set conformed by $\eta_1(x)=\widehat{\psi}_4(x)-\widehat{\zeta}_4(x)=\frac{3}{2}\widehat{\psi}_4(x)$ and $\eta_2(x)=\widehat{\psi}_3(x)-\widehat{\zeta}_3(x)=\frac{3}{2}\widehat{\psi}_3(x)$ in Eq. \eqref{eq:sed2o}. Now, the solutions of Eq. \eqref{eq:ggen} are given by the expression $$\xi_k(x)=\sqrt[4]{xC_1(x)}\eta_k(x)=\frac{3}{2}\sqrt[4]{xC_1(x)}\widehat{\psi}_{5-k}(x), \quad k=1,2,$$ which is given as follows:
\begin{equation}\label{solRx0}
\begin{array}{ccr}
\xi_1(x)&=&\sqrt[4]{\frac{C_1(x)}{x}}\cdot \sqrt{C_2(x)}  \cdot \exp\left(\mu\int\frac{\sqrt{x}dx}{C_2(x)\sqrt{C_1(x)}}\right)\\&&\\
\xi_2(x)&=&\sqrt[4]{\frac{C_1(x)}{x}}\cdot \sqrt{C_2(x)} \cdot \exp\left(-\mu\int\frac{\sqrt{x}dx}{C_2(x)\sqrt{C_1(x)}}\right)
\end{array}
\end{equation}
Now we are interested in the obtaining of the identity connected component of the Galois group for Eq. \eqref{eq:ggen} to apply Theorem \ref{tmora1}, see \ref{sec:DGT} for details and theoretical background. We recall that the differential field for Eq.  \eqref{eq:ggen} is $K=\mathbb{C}(x)$ and the corresponding Picard Vessiot extension is $L_1=K( \xi_1(x),\xi_1'(x))$. We observe that the logarithmic derivatives of the solutions are algebraic functions, therefore the differential Galois Group $G_1:=G(L_1/K)$ of $VE_1$ is not a connected group. We denote by $\overline{K}$ the algebraic functions over $K$, contained in the Picard-Vessiot $L_1$. Therefore, the differential Galois group over $L_1/\overline{K}$ is a connected group, i.e.,  $(G_1)^0=G(L_1/\overline{K})$. Now, to compute the Galois group $G(L_1/\overline{K})$ we start writing Eq. \eqref{eq:ggen} as a linear differential system, thus we obtain
\begin{equation}\label{eq:fvsys}
\begin{pmatrix}
\xi\\\xi'
\end{pmatrix}'=\begin{pmatrix}
0&1\\g(x)&0
\end{pmatrix}\begin{pmatrix}
\xi\\\xi'
\end{pmatrix}
\end{equation}
The fundamental matrix of Eq. \eqref{eq:fvsys} is
\begin{equation}\label{eq:fmfve}
    U=\begin{pmatrix}
    \xi_1&\xi_2\\
    \xi_1'&\xi_2'
    \end{pmatrix}
\end{equation}
We compute the action of the Galois Group $G(L_1/\overline{K})$ over the fundamental matrix \eqref{eq:fmfve}, therefore $\sigma$, element of $G(L_1/\overline{K})$, satisfies $$\sigma U=\sigma\begin{pmatrix}\xi_1&\xi_2\\\xi_1'&\xi_2'\end{pmatrix}=\begin{pmatrix}\xi_1&\xi_2\\\xi_1'&\xi_2'\end{pmatrix} A_\sigma,\quad A_\sigma=\begin{pmatrix}\delta&0\\0&\delta^{-1}\end{pmatrix}\in \mathbb{G}_m.$$ 
Since $(\mathbb{G}_m,\cdot)$ is isomorphic to the multiplicative group $(\mathbb{C}^*,\cdot)$, we conclude that $G(L_1/\overline{K})$ is an abelian group. We recall that a  Picard-Vessiot extension is purely transcendental if and only if its Galois group is connected, see \ref{sec:DGT}. For the interested reader, in \ref{sec:kovacic} we present in detailed way the solutions and Galois group $G_1=G(L_1/K)$, with $K=\mathbb{C}(x)$, by the direct application of Kovacic Algorithm.\\

As we are looking for obstructions to the integrability of the axisymmetric Bianchi IX cosmological model, taking into account Theorem \ref{tmorasi}, we must to consider the second variational equation over $\overline{K}$, because unfortunately the identity component of the Galois group of the first variational equation is abelian. In this way, we should compute the second variational equation $VE_2$. To obtain it $VE_2$, see \ref{sec:DGT} and references \cite{morasi,apt},  we make the change of variables 
\begin{equation}\label{chanvarAB}
(A,B)\mapsto (A_s+\epsilon C_A+\epsilon^2 D_A,B_s+\epsilon C_B+\epsilon^2 D_B),
\end{equation} where $(As,Bs)$ is a particular solution of Eq. \eqref{hebi}. Next step is to apply the change of variable given in Equation \eqref{chanvarAB} into Equation \eqref{hebi}. Using the invariant plane $\Gamma=\Gamma_E$, provided by $A=B, \dot A=\dot B$ and collecting the coefficients of $\epsilon^2$ we obtain the second variational equation related with Eq. \eqref{bnveg}. The second variational equation related to Eq. \eqref{FVEG} comes from the collecting of coefficient of $\epsilon^2$ after the restriction of $\Gamma=\Gamma_E$ in Equation \eqref{bnveg}, written in terms of $\hat \psi$ and $\hat \zeta$, i.e. we get the second variational equation  along the Taub solution. Recall that $\Gamma=\Gamma_E$ is provided by Equation \eqref{eqalpha} (some computations were made with Maple).\\
Now, in a similar way as for Equation \eqref{nveg}, we write one equation coming from the difference $\varepsilon=\hat\psi-\hat\zeta$, which now reads
\begin{equation}\label{eqsve}
    \begin{array}{c}
   \ddot \varepsilon+\left(\frac{\dot x}{x}\right)\dot \varepsilon+\left(\frac{17}{8x^2}-\frac{\dot x^2}{2x^2}-\frac{1}{2}\lambda\right)\varepsilon=f(\chi_1,\dot\chi_1,\ddot\chi_1,\psi_4,\dot\psi_4,x,\dot x), \\ f(\chi_1,\dot\chi_1,\ddot\chi_1,\psi_4,\dot\psi_4,x,\dot x)=
 \left(\left(\frac{2\lambda}{x}+\frac{\dot x^2}{x^3}-\frac{17}{4 x^3}\right)\psi_4+ \frac{\dot x}{x^2}\dot \psi_4\right)\chi_1 -\frac{4}{x}\chi_1\ddot \chi_1\\+\left(\frac{2\lambda}{x}+\frac{\dot x^2}{x^3}-\frac{2}{x^3}\right)\chi_1^2-\frac{1}{x}\dot\chi_1^2-(3\frac{ \dot x}{x^2}\psi_4+\frac{1}{x}\dot \psi_4)\dot \chi_1-\frac{4}{x}\psi_4\ddot \chi_1 -2\frac{ \dot x}{x^2}\chi_1\dot\chi_1,
    \end{array}
\end{equation}
where $x=x(t)$, $\dot x=\dot x(t)$ is the Taub solution in $\Gamma_E$, $\psi_4$ is one solution of equation \eqref{FVEG} and $\chi_1$ is one solution of equation \eqref{nveg}. 
Applying again the algebrization process through the Hamiltonian change of variable $x=x(t)$, such as for Equation \eqref{eq:ave}, we obtain the algebraic form for equation \eqref{eqsve} as follows.
\begin{equation}\label{eqasve}
\begin{array}{l}
\varepsilon''+p(x)\varepsilon'+q(x)\varepsilon=
\frac{\widetilde{f}}{x^3C_1(x)},\\\widetilde{f}=(28\lambda x^3-27x+2E)\eta_1^2-12x^2(4\lambda x^2- 1)\eta_1\eta_1'\\-8x^2C_1(x)\eta_1\eta_1''+xC_1(x)\widehat{\psi_4}'\eta_1 +(28\lambda x^3-54x+2E)\widehat{\psi}_4\eta_1\\-x(28\lambda x^3-9x+2E)\widehat{\psi}_4\eta_1'
-x^2C_1(x)\widehat{\psi}_4'\eta_1'-4x^2C_1(x)\widehat{\psi}_4'\eta_1'',
\end{array}
\end{equation}
where $p(x)$ and $q(x)$ are given in Eq. \eqref{eq2o},  $\widehat{\psi}_4$ is solution of Eq. \eqref{eq:ave}  (given explicitly in Eq. \eqref{eq:foave}) and $\eta_1$ is one solution of Eq. \eqref{eq2o}. Moreover, due to $\eta_1''=-p\eta_1'-q\eta$ and $\widehat{\psi}_4=\frac{2}{3}\eta_1$ we can simplify the expression for $\widetilde{f}$ as follows
\begin{equation*}
\begin{array}{l}
\widetilde{f}=-\frac{1}{3}(52\lambda x^3-459x+14E)\eta_1^2\\-\frac{2x}{3}(4\lambda x^4-32\lambda x^3-3x^2+(2E+12)x-4E)(\eta_1')^2\\
-\frac{2}{3}(28\lambda x^4-4\lambda x^3 -21 x^2+(14E+99)x-2E)\eta_1\cdot \eta_1'
\end{array}
 \end{equation*}
Now, writing $\eta_1'$ in terms of $\eta_1$ and $x$ we have
$$\eta_1'=\left(\frac{xC_1(x)C_2'(x)-C_1(x)C_2(x)+2\mu x^{3/2}C_1(x)}{2xC_1(x)C_2(x)}\right)\eta_1$$
and therefore $\widetilde{f}$ becomes $$\widetilde{f}=\frac{\mathcal{P}(\sqrt{x})}{C_2(x)^2}\eta_1^2=\frac{\mathcal{P}(\sqrt{x})}{C_2(x)^2\sqrt{xC_1(x)}}\xi_1^2,\quad \displaystyle{\mathcal{P}(\sqrt{x})=\sum_{k=0}^{18}a_k(\sqrt{x})^k},$$ where $\eta_1=\frac{\xi_1}{\sqrt[4]{xC_1(x)}}$ and $\xi_1$ is given in Eq. \eqref{solRx0}. The coefficients of the polynomial $\mathcal{P}$ are 

$$\begin{array}{l}
a_{18}=-\frac{31850}{3}\lambda^3,\,\,a_{17}=0,\,\, a_{16}=4900\lambda^3,\,\, a_{15}=-560\lambda^2\mu,\,\, a_{14}=\frac{104125}{2}\lambda^2,\\ a_{13}=- \frac{1960}{3}\lambda^2\mu,\,\, a_{12}=\lambda((-2835E + 96775)\lambda - \frac{8}{3}\mu^2),\,\, a_{11}=-2940\lambda\mu,\\a_{10}=1120\lambda^2E + \frac{64}{3}\mu^2\lambda + 1771350\lambda,\,\, a_{9}=-344\lambda\mu(E + \frac{2975}{172}),\\a_8=(21560E + 632100)\lambda + 2\mu^2,\,\, a_7=-\frac{364}{3}\lambda E\mu + 2520\mu,\\
a_6=-232\lambda E^2-\frac{4}{3}\mu^2 + 15400\lambda-8\mu^2 + 6460650,\,\, a_5=(-1632E - 12180)\mu,\\ a_4=64\lambda E^2 + \frac{8}{3}E\mu^2 + 231420 E + 1367100,\,\,
a_3=-32E\mu(E + 16),\\    a_2=1784 E^2 + 52080 E,\,\, a_1=-\frac{16}{3}E^2\mu,\,\, a_0=- \frac{16}{3}E^2(E - 93).
\end{array}$$

In a similar way as for equations \eqref{eq2o} and \eqref{eq:ggen}, we obtain the invariant form of equation \eqref{eqasve} through the change of dependent variable 
$$\varepsilon=\frac{\omega}{\sqrt[4]{xC_1(x)}}$$
as follows.
\begin{equation}\label{eq:ggensve}
\omega''=g(x)\omega+\frac{\mathcal{P}(\sqrt{x})}{x^{13/4}C_1(x)^{5/4}C_2(x)^2}\xi_1^2(x),
\end{equation} 
where $g(x)$ is given in equation \eqref{eq:ggen}.

We must solve an inhomogeneous second order differential equation, Eq. \eqref{eq:ggensve}.
Thus, computing the wronskian 
\begin{equation}\label{eq:wrons}
W:=W(\xi_1(x),\xi_2(x))=18i\sqrt{6E^2\lambda+14700}
\end{equation} and using the method of variation of parameters we obtain as solutions of \eqref{eq:ggensve} the following
\begin{equation}\label{solvarpam}
\begin{array}{lll}
\omega(x)&=&c_1\xi_1(x)+c_2\xi_2(x)+\omega_p(x)\\\omega_p(x)&=&
-\frac{\xi_1(x)}{W}\int\frac{\mathcal{P}(\sqrt{x})}{x^{4}C_1(x)^{3/4}C_2(x)}\xi_1(x)dx+\\ && \frac{\xi_2(x)}{W}\int\frac{\mathcal{P}(\sqrt{x})}{x^{7/2}C_1(x)^{5/4}C_2(x)^2}\xi_1^3(x)dx
\end{array}
\end{equation}
It is convenient to write the inhomogeneous equation, Eq. \eqref{eq:ggensve}, as an homogeneous one, using as a new variable $\varphi=\xi_1^2$. Thus,
\begin{equation}\label{eq:omegassp}
\begin{array}{lll}
\omega''&=&g(x)\omega+s(x)\varphi(x),\quad s(x)= \frac{\mathcal{P}(\sqrt{x})}{x^{7/2}C_1(x)^{5/4}C_2(x)^2}\\
\varphi'&=&2h(x)\varphi,\quad h(x)=(\log(\xi_1(x)))'
\end{array}    
\end{equation}
The key point is $h(x)=(\log(\xi_1(x)))'$ belongs to $\overline{K}=\overline{\mathbb{C}(x)}$, i.e., the field of coefficients of Eq. \eqref{eq:omegassp} becomes $\overline{K}$. 
Rewritten it in matrix form
\begin{equation}\label{eq:matssw}
\begin{pmatrix}\omega\\ \omega'\\ \varphi \end{pmatrix}'=\begin{pmatrix}0&1&0\\g(x)&0&s(x)\\ 0&0&2h(x)\end{pmatrix}\begin{pmatrix}\omega\\ \omega'\\\varphi \end{pmatrix}.
\end{equation}
The fundamental matrix of Eq. \eqref{eq:matssw} is
\begin{equation}\label{eq:matfund}
\Phi=\begin{pmatrix}\xi_1&\xi_2&\omega_p\\\xi_1'&\xi_2'&\omega_p'\\ 0&0&\varphi \end{pmatrix}.
\end{equation}
The Picard-Vessiot Extension of Eq. \eqref{eq:matssw} is $\overline{K}(\xi_1,\xi_1',\omega_p,\omega_p')\supset \overline{K}$. We compute the action of the Galois Group $G$ over the fundamental matrix \eqref{eq:matfund}, recall that $\sigma$, element of $G$, satisfies $\sigma(\xi_1)=\delta{\xi_1}$ and $\sigma(\xi_2)=\delta^{-1}\xi_2$. Therefore $\sigma(\varphi)=\sigma(\xi_1^2)=\delta^2\varphi$. Finally, setting $$I_1:=\int\frac{\mathcal{P}(\sqrt{x})}{x^{4}C_1(x)^{3/4}C_2(x)}\xi_1(x)dx:=\int y_1(x)\xi_1(x)dx,$$
$$I_2:=\int\frac{\mathcal{P}(\sqrt{x})}{x^{7/2}C_1(x)^{5/4}C_2(x)^2}\xi_1^3(x)dx:=\int y_2(x)\xi_1^3(x)dx,$$ we observe $\sigma(I_1')=y_1\sigma(\xi_1)$ and $\sigma(I_2')=y_2\sigma(\xi_1^3)$, because $y_i\in \overline{K}$, $i=1,2$. Hence, $\sigma(I_1)=\delta I_1+\gamma_1$ and $\sigma(I_2)=\delta^3I_2+\gamma_2$. Taking into account the expression for $\omega_p$ in Eq. \eqref{solvarpam}, as $W \in K\subset \overline{K}$, we obtain
$$\sigma(\omega_p)=\sigma\left(-\frac{\xi_1}{W}I_1+\frac{\xi_1^3}{W}I_2\right)=\delta\gamma_1\xi_1+\delta^{-1}\gamma_2\xi_2+\delta^2\omega_p,$$ where $\delta\in \mathbb{C}^*,\,\gamma_1,\,\gamma_2\in \mathbb{C}$.
Therefore, $\sigma$ satisfies $$\sigma\begin{pmatrix}\xi_1&\xi_2&\omega_p\\\xi_1'&\xi_2'&\omega_p'\\ 0&0&\varphi \end{pmatrix}=\begin{pmatrix}\xi_1&\xi_2&\omega_p\\\xi_1'&\xi_2'&\omega_p'\\ 0&0&\varphi \end{pmatrix} B_\sigma,\quad B_\sigma=\begin{pmatrix}\delta&0&\delta\gamma_1\\0&\delta^{-1}&\delta^{-1}\gamma_2\\ 0&0&\delta^2 \end{pmatrix}.$$ 
The matrix $B_\sigma\in G$ can be written as follows
$$B_\sigma=\begin{pmatrix}\delta&0&0\\0&\delta^{-1}&0\\ 0&0&\delta^2 \end{pmatrix}\cdot  \begin{pmatrix}1&0&\gamma_1\\0&1&\gamma_2\\ 0&0&1 \end{pmatrix}.$$
According to \cite[\S 8.4]{humpreys}, $G = D \ltimes N$, is the semidirect product between the Diagonal subgroup of $G$ (left factor in the decomposition of $B_\sigma)$ and the Normal Unipotent subgroup of $G$ (right factor in the decomposition of $B_\sigma$).

Then we state the following conjecture.\medskip

\noindent\textbf{Conjecture.} \textit{If generically $\lambda\neq 0$ and $E\neq 0$, then $\gamma_1$ and $\gamma_2$ are not simultaneously zero.}\medskip

By the differential Galois theory the above conjecture means that either $I_1$ or $I_2$ do not belong to $L_1=\mathbb{C}(z,\xi_1,\xi_1')$, i.e., these integrals are not expressed as algebraic functions of $z$, $\xi_1$ and $\xi_1'$. This conjecture is natural because, as usually happens, when you integrate functions in a differential field you do not obtain functions in the same field. Moreover, the proof of this conjecture seems to us that it would be quite involved: as far as we know there are not algorithms, in the literature, to solve this problem. More elementary situations were studied in the papers \cite{almp2,ay}.\\

As the Galois group $G$ is connected and non abelian, using Theorem \ref{tmorasi}, we obtain the result that we are looking for: 

\begin{theorem}\label{mainth}
Assuming the above conjecture, the Axisymmetric Bianchi IX cosmological model is not meromorphically integrable for generic values of $(E,\lambda)$.
\end{theorem}

\section{Numerical experiments}

To illustrate the conjecture about the non integrability of  the axisymmetric Bianchi IX models with cosmological constant and reduced to two
degrees of freedom, we calculate Poincar\'e sections through the surface $P_B=0$ and $\dot {P_B}>0$.
It is not our intention to make a through investigation of the section phase space,  so we have
chosen to view selected values of  the  $(\lambda, energy)$ pair.

\begin{figure}
\centering
\begin{tabular}{l}
\includegraphics[width=6cm]{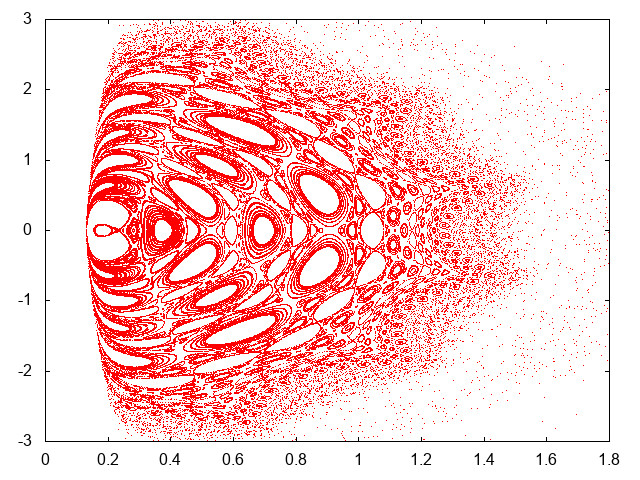}
\includegraphics[width=6cm]{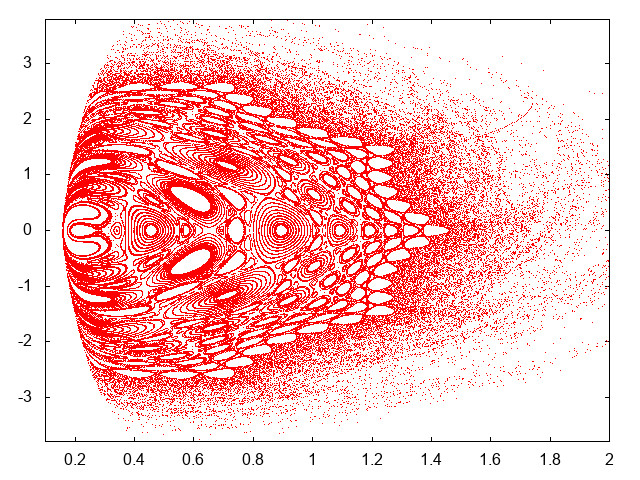}\\
\includegraphics[width=6cm]{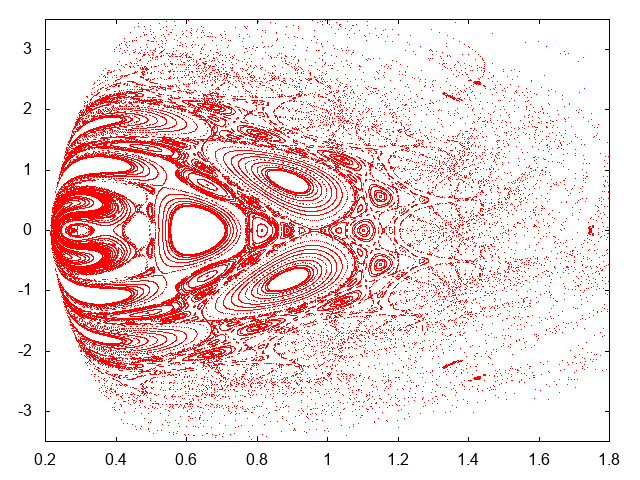}
\includegraphics[width=6cm]{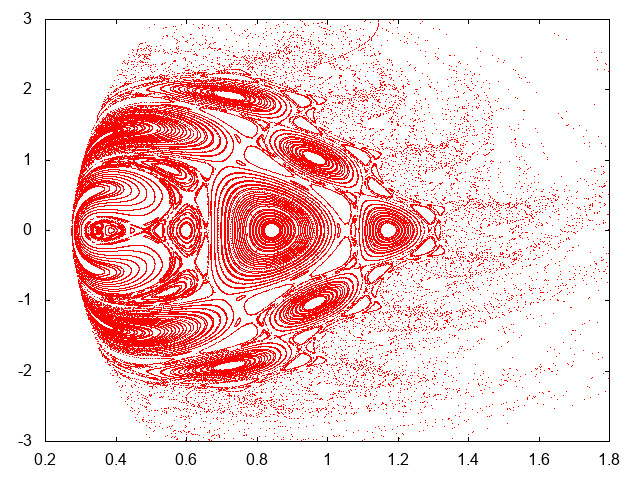}
\end{tabular}
\caption{ $(A,P_A)$ projection of the section $P_B=0$, $\dot P_B>0$. $\lambda$ is fixed at $0.1$ and the energy has the following values $0.25$, $0.3$ (top) and $0.4$, $0.5$ (bottom).
 Note the increase of the island region as the energy increases.}
\label{varenergy}
\end{figure}

First of all let us look at the behaviour of the dynamics taking a fixed value of  $\lambda$, $\lambda=0.1$,
and varying the  energy: $0.25,~ 0.3,~ 0.4~and ~0.5$. We have left in all figures some island regions
 empty to make the viewing better.  Around $A=1$ the iterates scape to infinity to the  De Sitter attractor.
 When the value of the energy increases the islands become bigger. This feature  appears more clearly
 in the bottom  panels of  fig  \ref {varenergy} which have  higher energy values

\begin{figure}
\centering
\begin{tabular}{l}
\includegraphics[width=6cm]{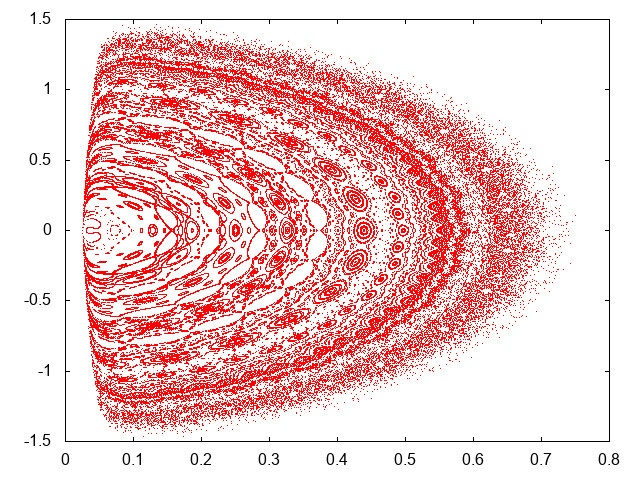}
\includegraphics[width=6cm]{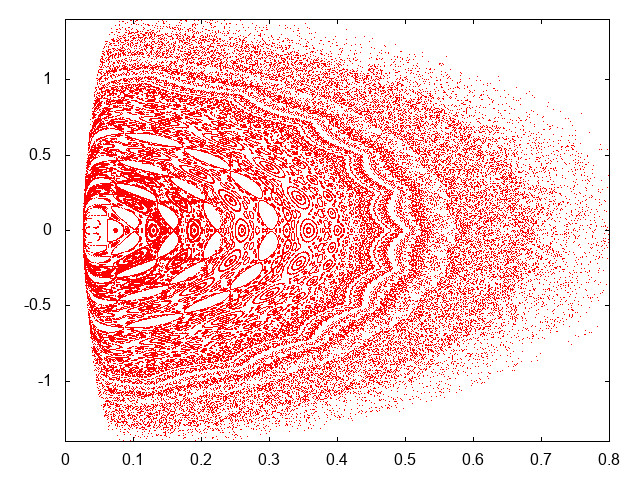}\\
\includegraphics[width=6cm]{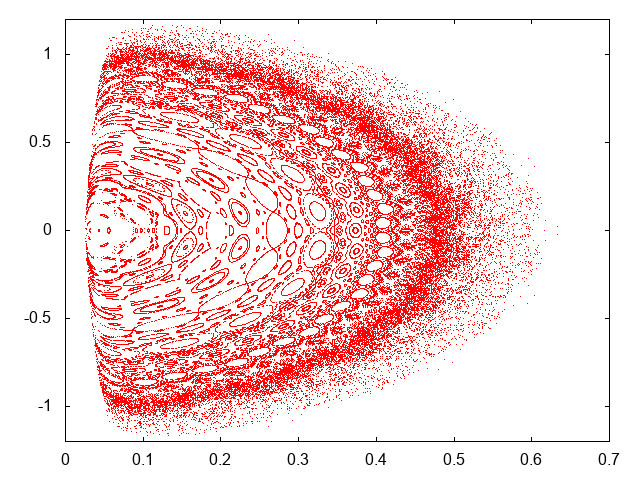}
\includegraphics[width=6cm]{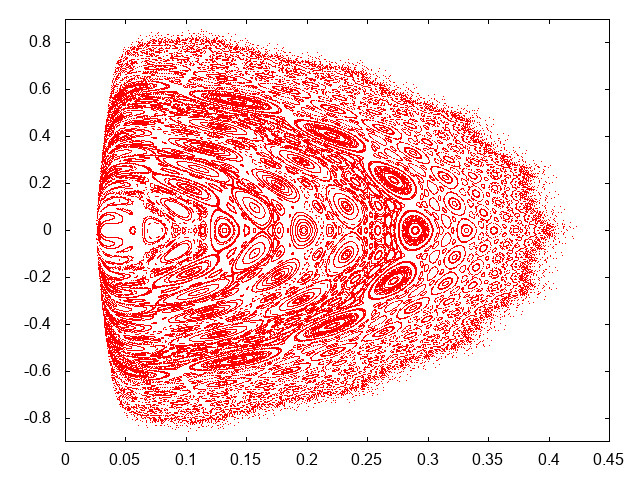}
\end{tabular}
\caption{ ($A,P_A$) projection. The energy is  fixed at $0.05$. Top, from left to right : $\lambda=0.6$ and $0.8$ ; bottom:  $\lambda=1.0$ and $1.2$ from left to right.
The region with more islands is the one with larger value of $\lambda$.}
\label{lamvar}
\end{figure}

One could want to know how the sections vary with a fixed value of the energy and different values of $\lambda$.
In fig  \ref {lamvar}  top  we have chosen $E=0.05$ and varied $\lambda$  from $0.6$ 
to $1.2$.  Note that the quantity of islands increase as $\lambda$ increases. As in fig \ref {varenergy} the border of
the Poincar\'e section are  escaping points  to larger regions of the phase space.

\begin{figure}
\centering
\begin{tabular}{l}
\includegraphics[width=6cm]{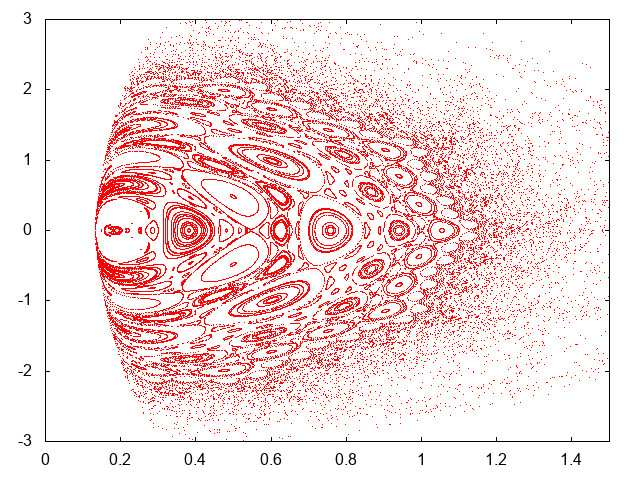}
\includegraphics[width=6cm]{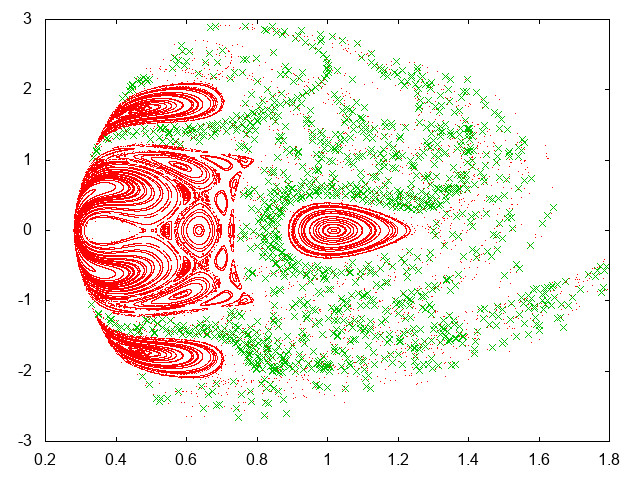}\\
\includegraphics[width=6cm]{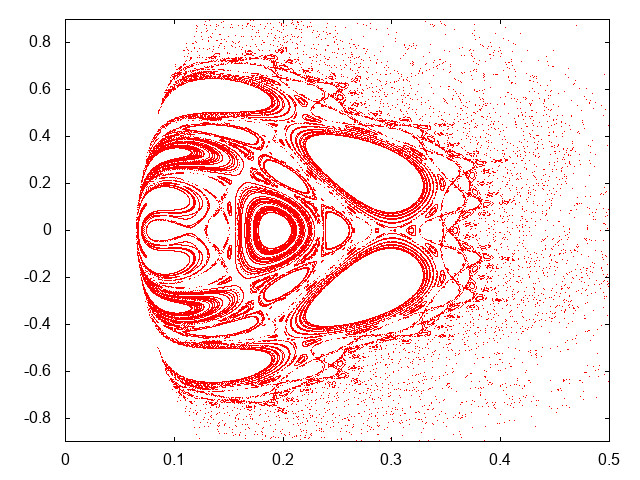}
\includegraphics[width=6cm]{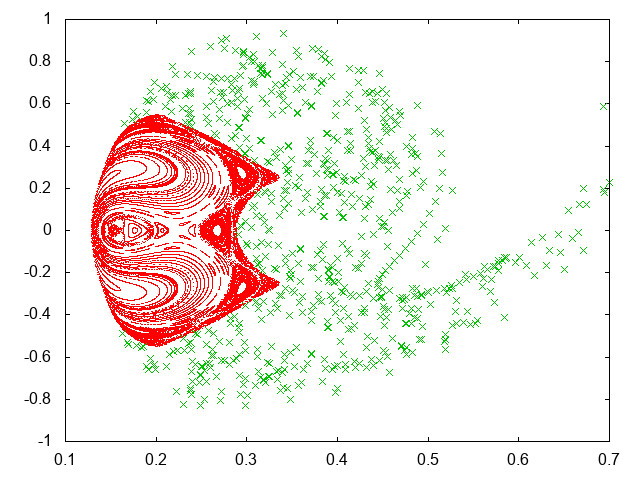}
\end{tabular}
{\caption{ Top: $\lambda=1.5$  and $energy=0.25$ (left) and $0.5$ (right). Bottom  $\lambda=1.2$
with $energy=0.12061132$ (left) and $0.22061113$ (right). The more organized region is destroyed and there
tongs of iterates embracing the reduced island region.}}
\label{escape}
\end{figure}

In fig \ref{escape} we show the sections when most of the iterates have escaped and return to the still ordered region;
 in fig \ref{escape} top  $\lambda=0.15$ and the energy varies from $0.25$ to $0.5$; in the bottom case $\lambda=1.2$ and
 the energy goes from $0.12061132$ to $0.22061132$. Note that the iterates which are not in the islands have the appearance
 of stable or unstable manifold s embracing the island region. This is the point to  recall that there is a saddle centre
  equilibrium point in the dynamics in $A=B=\frac{1}{\sqrt{ 4\lambda}}=E$.  This point is a scatterer for
  the dynamics taking away and bringing back iterates of the Poincar\'e section.

\begin{figure}
\centering
\begin{tabular}{l}
\includegraphics[width=6cm]{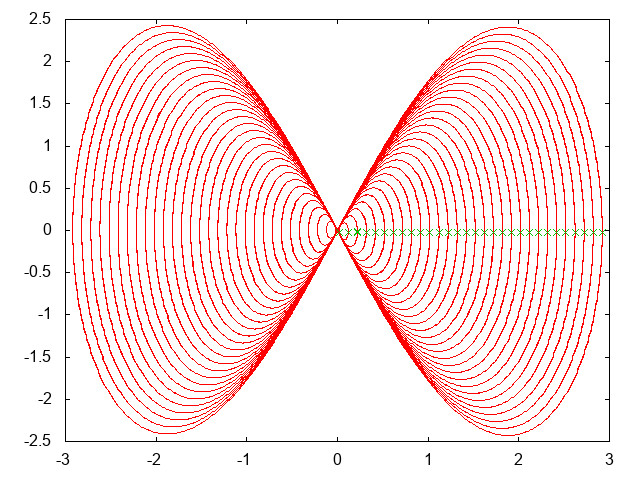}
\includegraphics[width=6cm]{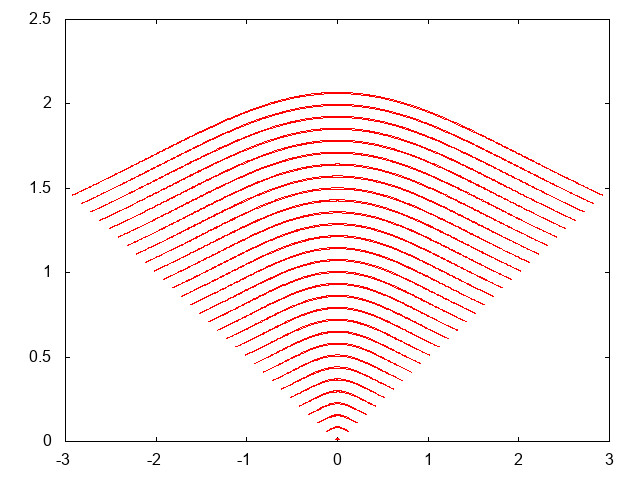}\\
\end{tabular}
\caption{ Right:$(A,P_A)$  and left $(A,B)$  projection. Compare the scales with these of the previous figure. Since the critical point was moved to infinity by taking $\lambda=0$ the family of
periodic orbits remains all over the phase space. Obviously there is no escaping orbits.}
\label{E0l0}
\end{figure}

An interesting feature is found  for the system when $\lambda=0$. The fixed point is in this case $A=B=\frac{1}{\sqrt{ 4\lambda}}=E$ and $E\rightarrow \infty$ when $\lambda\rightarrow 0$. In  Fig \ref{E0l0} we  show 30000 initial conditions integrated each 300.  The behaviour is continuos and no trajectory scape to infinity.

\begin{figure}
\centering
\begin{tabular}{l}
\includegraphics[width=6cm]{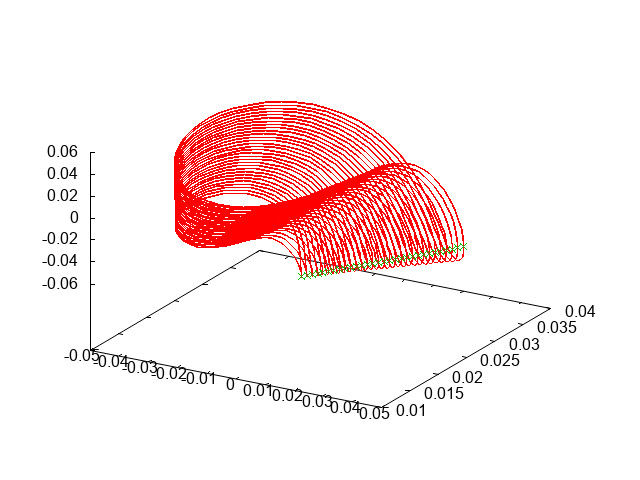}
\includegraphics[width=6cm]{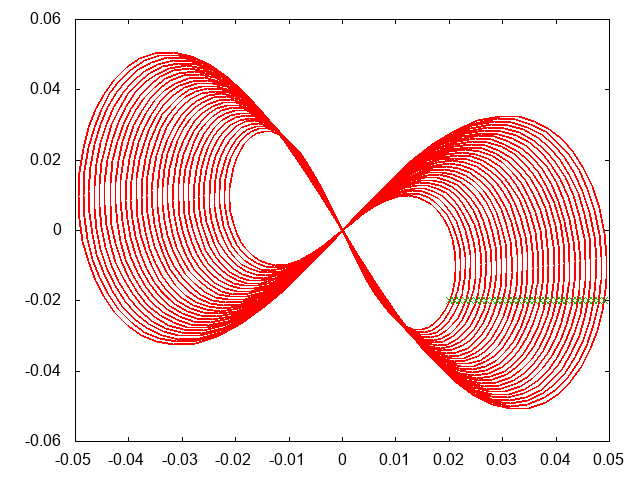}\\
\includegraphics[width=6cm]{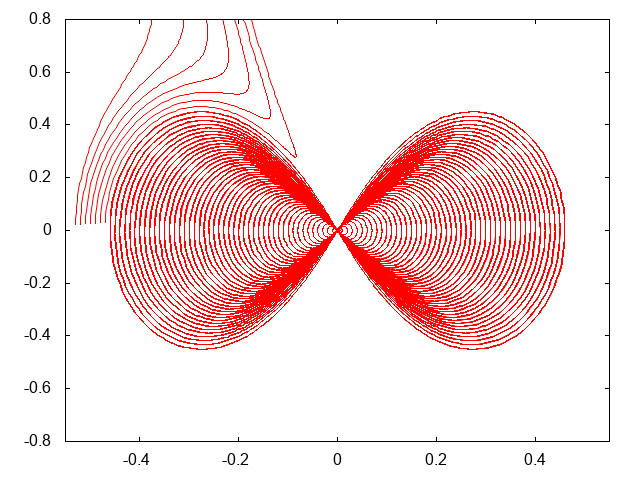}
\includegraphics[width=6cm]{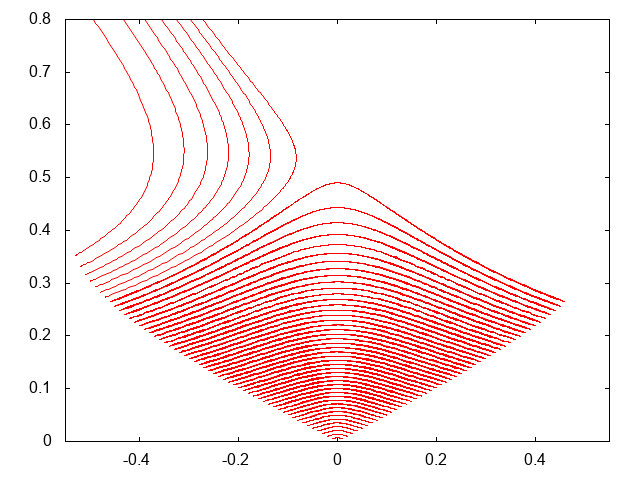}
\end{tabular}
\caption{ Top: left, $ABP_A$ projection of some trajectories with zero energy and $\lambda= 3.5$; right, the same apparent periodic trajectories projected
 in the plane $(A,P_A)$. Bottom: left, the periodic orbits are interrupted by escaping trajectories around $0.5$. The saddle centre
is located at $A_f=-0.53452247$; right, the $(A,B)$ projection showing the exit in the scale $B$.}
\label{periof}
\end{figure}
 
\appendix	

\section{Differential Galois Theory and Hamiltonians}
\label{sec:DGT}

Differential Galois Theory, see \cite{crha,vasi} among others, also known as Picard-Vessiot theory, is the Galois theory of linear differential
equations. We recall that rational functions over $\mathbb{C}$ are denoted by $\mathbb{C}(x)$. Consider differential equations in the form $\mathcal L:= y''+ay'+by=0$, where $a,b\in K:=\mathbb{C}(x)$, usually called the differential coefficient field of  $\mathcal L$. Suppose that $y_1, y_2$ is a basis of solutions of $\mathcal L$. Let $F=K(y_1, y_2, y_1', y_2')$ be the differential
	extension of $K$ such that $\mathbb{C}$ is the field of constants
	for $K$ and $F$. In these terms, we say that $F$, the smallest
	differential field containing $K$ and $\{y_{1},y_{2}, y_1',y_2'\}$, is the
	\textit{Picard-Vessiot extension} of $K$ for $\mathcal L$. Therefore, the
	group of all $K$-differential automorphisms of $F$ over $K$ is called the
	{\it{differential Galois group}}, or simply the \textit{Galois group}, of $F$ over $K$ and is denoted by
	$G(F/K)$. In terms of differential Galois theory the notion of integrability is the following:

 We say that the differential equation $\mathcal L$ is \textit{integrable} if and only if the Picard-Vessiot
	extension $F\supset K$ is obtained as a tower of differential
	fields $K=F_0\subset F_1\subset\cdots\subset F_m=F$ such that
	$F_i=F_{i-1}(\eta)$ for $i=1,\ldots,m$, where either
	\begin{enumerate}
		\item $\eta$ is {\emph{algebraic}} over $F_{i-1}$, that is $\eta$ satisfies a
		polynomial equation with coefficients in $F_{i-1}$.
		\item $\eta$ is {\emph{primitive}} over $F_{i-1}$, that is $\eta' \in F_{i-1}$.
		\item $\eta$ is {\emph{exponential}} over $F_{i-1}$, that is $\eta' /\eta \in F_{i-1}$.
	\end{enumerate}
Kolchin Theorem states the relation between the integrability of of $\mathcal L$ and the \emph{solvability} of the connected identity component of its differential Galois group as follows:\\

\noindent \textbf{Kolchin Theorem.}\emph{
	The equation $\mathcal L$ is integrable if
	and only if the connected identity component of the differential Galois group is a solvable group.}\\

We denote by $(G(F/K))^0$ the connected identity component of the differential Galois group of the extension $F$ over the differential field $K$. Picard-Vessiot extension is purely transcendental if and only if its Galois group is connected, that is, $(G(F/K))^0=G(F/K)$, see \cite{crha} and references therein for further details.\\

The reduced form, also known as invariant normal form, of the equation $\mathcal L$ is as follows:

$$\mathcal R:=\zeta''=r\zeta,\quad r\in K.$$

We recall that equation $\mathcal R$ can be obtained from equation
$\mathcal L$ through the change of variable
$$y=e^{-\frac{1}{2}\int_{}^{}a}\zeta,\quad r=\frac{a^2}{4}+\frac{a'}{2}-b.$$

Through the change of variable
$v=\zeta'/\zeta$ we get the associated Riccati equation
to equation $\mathcal R$ as follows:
$$v'=r-v^2,\quad v=\frac{\zeta'}{\zeta},$$
where $r$ is given by the above equation. Moreover, the above Riccatti equation has one algebraic solution
over the differential field $K$ if and only if the differential
equation $\mathcal R$ is integrable.

The concepts related with integrability of Hamiltonian Systems are very well known. The Hamilton Equations for the Hamiltonian $H$ are 

$$
\dot x_i=\partial H/\partial y_i,
\dot y_i=-\partial
H/\partial x_i, \, i=1,...,n
$$

 One says that $X_H=(\partial H/\partial y_i,-\partial
H/\partial x_i)$, $i=1,...,n$,  is completely integrable or
Liouville integrable if there are $n$ functions $f_1=H$,
$f_2$,..., $f_n$, such that

\vskip 0.2cm

(1) they are functionally independent i.e., the 1-forms $df_i$
$i=1,2,...,n$, are linearly independent over a dense open set
$U\subset  M$, $\bar  U= M$;

\vskip 0.15cm

(2) they form an involutive set, $\{f_i,f_j\}=0$, $i,\,
j=1,2,...,n$.

\vskip 0.2cm

We remark that in virtue of item (2) above the functions $f_i$,
$i=1,...,n$ are first integrals of $X_H$. It is very important to
be precise regarding the degree of regularity of these first
integrals. Here we assume that the first integrals
are meromorphic. 

Given a dynamical system of first order,

\begin{equation*}\label{sd}
\dot z= X_H(z),
\end{equation*}
we can write the variational equations  along a particular
integral curve  $z=\phi(t) $  of the Hamiltonian vector field $X_H$

\begin{equation*}
\dot\xi= X_H^\prime(\phi(t))\xi. \label{ve}
\end{equation*}
We denote by $\Gamma$ the Riemann surface representing the integral curve $z=\phi(t)$ which is not an equilibrium point of the vector field $X_H$. The following theorem is included in a new Theory of Integrability of Dynamical Systems using Differential Galois Theory, see  \cite {mora1, mora2} and  see also
	\cite{mo}:\\

\begin{theorem}[Meromorphically non-integrability theorem, \cite{mora1}]\label{tmora1}
Assume  a complex analytic
	Hamiltonian system is meromorphically completely integrable  in a
	neighborhood of the integral curve $z=\phi(t)$ . Then the identity
	component of the Galois group of the variational equation is a commutative group.
	\end{theorem}

The previous theorem has been generalized  to higher order
variational equations
$\mathrm{VE}_k$ along
$\Gamma$, with $k>1$, $VE_1$ being the first variational equation VE.\medskip

The ``fundamental'' solution of $\mathrm{VE}_k$ of our dynamical system is given by
$$(\phi^{(1)}(t),\phi^{(2)}(t),\ldots, \phi^{(k)}(t)),$$ where
$$\phi(z,t)=\phi(z_0,t)+\phi^{(1)}(t)(z-z_0)+\ldots+ \displaystyle{\frac{1}{k!}}\phi^{(k)}(t)
(z-z_0)^k+\ldots$$ the Taylor series up to order $k$ of the flow
$\phi(z,t)$ with respect to the variable $z$ at the point
$(z_0,t)$. That is, $\phi^{(k)}(t)=\displaystyle{\frac{\partial^k}
	{\partial z^k}}\phi(z_0,t).$ The initial conditions are
$\phi^{(1)}(0)=Id_m$ and $\phi^{(j)}(0)=0$ for all $j>1.$ 

Since through a linearization process   we can  consider the equations $\mathrm{VE}_k$ as linear
differential equations,  we can talk about their Picard-Vessiot
extensions and about their Galois groups $G_k$. Thus, we arrive to the following theorem, see \cite{morasi}.

\begin{theorem}[Higher order variational equation theorem, \cite{morasi}]\label{tmorasi}
Assume that a complex analytical Hamiltonian system
	is   integrable  by meromorphic first integrals in a neighborhood
	of the integral curve $z=\phi(t)$. Then
	the identity components $(G_k)^0$, $k\geq 1$,  of the
	Galois groups of the variational equations along $\overline\Gamma $ are
	commutative.
	\end{theorem}
	
Recently, in \cite{almp2,apt}, the higher order variational equation theorem where applied to obtain non-integrability of planar vector fields and Painlev\'e Equations. The explicit formulas presented in \cite{almp2,apt} are also useful here. In particular, to obtain the first and second variational equation of a Hamiltonian system $\dot z=X_H$, with particular solution $\gamma(t)$, we make the changes of variables $$x_i=x_{is}(t)+\epsilon w_i+\epsilon^2 r_i,\, y_i=y_{is}(t)+\epsilon \chi_i+\epsilon^2 \phi_i,\, \gamma(t)=(x_{1s},\ldots,x_{ns},y_{1s},\ldots,y_{ns}),$$ into the Hamilton equations $\dot z=X_H(z)$. The first variational equation is obtained collecting the coefficients of $\epsilon$ in the previous Hamilton equations, while the second variational equation is obtained by the collecting of the coefficient of $\epsilon^2$ in the previous Hamilton equations.
\section{Application of Kovacic Algorithm}
\label{sec:kovacic}
In order to apply Kovacic Algorithm, see \cite{kovacic:1986} and \cite{acbl},  for the general case $\Delta\neq 0$, $E\neq 0$ and $\lambda\neq 0$, we start considering the set of poles of $g$ that including the point at infinity is given by $\Upsilon= \{0,\rho_0,\rho_1,\rho_2,\infty\}$. Moreover, we observe that $\circ(g_c)=2$ for all  $c\in\Upsilon$, i.e., the order of each pole of $g$ is $2$. Thus, we could fall in any case of Kovacic Algorithm and we should check conditions $c_2$ and $\infty_2$ for each of them. That is, Laurent expansion of $g$ around each $c\in\Upsilon$. Therefore we have:
$$g(x)=\frac{b_0}{x^2}+\ldots,\quad \textrm{around } x=0,\quad b_0=\frac{5}{16}$$
$$g(x)=b_\infty{x^2}+\ldots,\quad \textrm{around } x=\infty,\quad b_\infty=2$$
$$g(x)=\frac{b_{\rho_1}}{(x-\rho_1)^2}+\ldots,\quad \textrm{around } x=\rho_1,\quad b_{\rho_1}\textrm{ is given by}$$
$$\frac{\lambda^4\kappa^2(128\lambda^2\kappa^4\rho_1^6-552\lambda\rho_1^4+128E\lambda\rho_1^3+315\rho_1^2-222E\rho_1+5E^2)}{(3(\kappa^{2/3}+\lambda)(\kappa^{4/3}+\lambda\kappa^{2/3}+\lambda^2))^2}.$$
In a similar way we obtain expressions for the coefficients of $x^{-2}$ in the Laurent expansion of $g$ around $\rho_2$ and $\rho_3$ respectively. Thus, $b_{\rho_2}$ and $b_{\rho_3}$ will include  algebraic combinations of $\rho_1$, $\rho_2$ and $\rho_3$ with some terms of $b_{\rho_1}$.\\
For suitability, we start with case 1 and after we follow with case 3 to discard them, showing that solutions are provided by case 2. Thus, starting with case 1 of Kovacic we obtain $[\sqrt{g}]_c=0$ and $\alpha^{\pm}_c=\frac{1\pm\sqrt{1+4b_c}}{2}$. In particular, $$\alpha^{+}_0=\frac{5}{4},\quad \alpha^{-}_0=-\frac{1}{4},\quad \alpha^{+}_\infty=2,\quad \alpha^{-}_\infty=-1,\quad \alpha^{\pm}_{\rho_j}\notin \mathbb{Q}.$$ By step 2 we obtain that does not exists a polynomial of degree $d$ satisfying at least one of the following options:
$$d=\alpha_\infty^\pm-(\alpha_{0}^\pm+\alpha_{\rho_1}^\pm+\alpha_{\rho_2}^\pm+\alpha_{\rho_3}^\pm)\in\mathbb{Z}_+.$$
Therefore, case 1 of Kovacic Algorithm is discarded and now we consider the case 3. By step 1 we have the same values for $b_c$ as for the previous cases. Thus, for $c\in \Upsilon\setminus\{\infty\}$ and $n\in\{4,6,12\}$ we have $$E_c=\{6\pm k\sqrt{1+4b_c}:\, 0\leq k\leq 6\}\cap\mathbb{Z},\, E_\infty^{(n)}=\left\{6\pm \frac{12k}{n}\sqrt{1+4b_\infty}:\, 0\leq k\leq 6,\right\}\cap \mathbb{Z}$$ Therefore we obtain the following sets

$$E_0=\left\{3k-3:\, 0\leq k\leq 6\right\},\,\, E_{\rho_j}=\{6\},\,\,E_\infty^{(n)}=\left\{6\pm \frac{36k}{n}:\, 0\leq k\leq 6\right\},\,\, j=1,2,3.$$ In particular we have $$E_\infty^{(4)}=\left\{6\pm 9k:\, 0\leq k\leq 6\right\},\quad E_\infty^{(6)}=\left\{6\pm 6k:\, 0\leq k\leq 6\right\},\quad E_\infty^{(12)}=\left\{6\pm 3k:\, 0\leq k\leq 6\right\}.$$
By step two, the set $D\neq \emptyset$ is as follows:
$$D=\left\{d\in\mathbb{Z}_+:d=\frac{n}{12}(e^{(n)}_{\infty}-e_0-e_{\rho_1}-e_{\rho_2}-e_{\rho_3})\right\}.$$ That is, there exist some non-negative integer values of $d$ satisfying
$$d=3k-\frac{n}{4}(k+3),\quad 0\leq k\leq 6.$$
Such values of $d$, depending of $k$ as in $E_c$ and in $E_\infty^{(n)}$, for $n\in\{4,6,12\}$  correspond to the following.
\begin{description}
    \item[$n=4$.] In this case, $d\in \{2k-3:\,\, 2\leq k\leq 6\}$, $e_0\in\left\{3k-3:\, 2\leq k\leq 6\right\}$, $e_\infty\in\left\{9k+6:\, 2\leq k\leq 6\right\}$.
    \item[$n=6$.] In this case, $d\in \{3k-9:\,\, 3\leq k\leq 4\}$, $e_0\in\left\{3k-3:\, 3\leq k\leq 4\right\}$, $e_\infty\in\left\{6k+6:\, 3\leq k\leq 4\right\}$,
    \item[$n=12$.] In this case we cannot obtain there are not values for $d\in \mathbb{Z}_+$.
\end{description}. By step 2 we obtain the polynomial $S=x(x-\rho_{1})(x-\rho_2)(x-\rho_3)$ and the rational function $\theta$ is given by $$\theta=\frac{n}{12}\sum_{c\in\Upsilon\setminus\{\infty\}}\frac{e_c}{x-c}=\frac{n}{12}\left(\frac{e_0}{x}+\frac{e_{\rho_1}}{x-\rho_1}+\frac{e_{\rho_2}}{x-\rho_2}+\frac{e_{\rho_3}}{x-\rho_3}\right).$$ Finally, by step 3, we set $P_n=-P_d$ and $$P_{i-1}=-SP'_i +((n-i)S'-S\theta)P_i-(n-i)(i+1)S^2gP_{i+1}.$$ Due to this recurrence is not satisfied, because we arrive to $P_{-1}\neq 0$ for all $d$ obtained previously with $n\in\{4,6\}$, we discard case 3 of Kovacic Algorithm. Due to case 1 and case 3 were discarded, we proceed to analyse case 2 of Kovacic Algorithm.\\

By step 1 of case 2 we obtain $$E_c=\{2,2-\sqrt{1+4b_c},2+\sqrt{1+4b_c},\,\forall c\in \Upsilon\}\cap \mathbb{Z}.$$ Now we should find the set $D$ corresponding to the possible degrees of the polynomials $P_d$. Then $$D=\left\{d\in\mathbb{Z}_+:\, d=\frac{1}{2}(e_{\infty}-e_0-e_{\rho_1}-e_{\rho_2}-e_{\rho_3})\right\},$$ $$e_0\in\{-1,2,5\},\quad e_\infty\in\{-4,2,8\},\quad e_{\rho_j}=2.$$ Therefore $D=\{0\}$, because the only one option to have a non-negative integer for $d$ is setting $e_\infty=8$ and $e_c=2$ for all $c\in\Upsilon\setminus \{\infty\}$. For instance the polynomial $P_0(x)=1$ and the rational function $\theta$, corresponding to the step 2 is given by $$\theta=\sum_{c\in\Upsilon\setminus\{\infty\}}\frac{1}{2}\left(\frac{e_c}{x-c}\right)=\frac{1}{x}+\frac{1}{x-\rho_1}+\frac{1}{x-\rho_2}+\frac{1}{x-\rho_3}.$$ Finally, by step 3 we have that $P_0$ and $\theta$ satisfy $$\theta''+3\theta\theta'+\theta^3-4g\theta-2g'\equiv 0.$$ Thus, the solution of the solutions $\xi_1$ and $\xi_2$ given by $$\xi_{1,2}=e^{\int \omega},\quad \mathrm{where}\,\, \omega^2-\theta\omega+\frac{1}{2}(\omega'+\omega^2-2g)=0.$$ For instance, the explicit solutions are given by
\begin{equation}\label{solRxk}
\begin{array}{ccr}
\xi_1&=&\sqrt[4]{\frac{C_1(x)}{x}}\cdot \sqrt{C_2(x)} \\ &&\\ &&  \cdot \exp\left(i\int {\sqrt{\frac{486x(E^2\lambda+2450)}{C_1(x)C_2(x)^2}}dx}\right)\\&&\\
\xi_2&=&\sqrt[4]{\frac{C_1(x)}{x}}\cdot \sqrt{C_2(x)} \\ &&\\ &&  \cdot \exp\left(-i\int {\sqrt{\frac{486x(E^2\lambda+2450)}{C_1(x)C_2(x)^2}}dx}\right)
\end{array}
\end{equation}
We recall that the differential field for equation \eqref{eq:ggen} is $K=\mathbb{C}(x)$. Thus, the corresponding Picard Vessiot extension is $L=\mathbb{C}(x)\langle \xi_1(x),\xi_1'(x)\rangle$. Computing the differential Galois group, we observe that for $\sigma\in Gal(L/\mathbb{C}(x))$, it satisfies $$\sigma^2(\xi_1(x)\cdot \xi_2(x))=(\xi_1(x)\cdot \xi_2(x))^2,\quad  \sigma\left(\left(\log\left(\xi_1(x)\xi_2(x)\right)\right)'\right)=\left(\log\left(\xi_1(x)\xi_2(x)\right)\right)',$$ $$\sigma^2\left(\left(\log\left(\frac{\xi_2(x)}{\xi_1(x)}\right)\right)'\right)=\left(\left(\log\left(\frac{\xi_2(x)}{\xi_1(x)}\right)\right)'\right)^2,$$ $$\xi_1(x)\cdot \xi_2(x)\notin\mathbb{C}(x),\quad\left(\log\left(\frac{\xi_2(x)}{\xi_1(x)}\right)\right)'\notin \mathbb{C}(x).$$ That is,
$$(\xi_1(x)\cdot \xi_2(x))^2=\frac{C_1(x)C_2(x)^2}{x}.$$
$$\left(\log\left(\xi_1(x)\cdot\xi_2(x)\right)\right)'=\frac{560\lambda^2 x^6+1365 \lambda x^4+191E\lambda x^3-630x^2+210Ex-4E^2}{xC_1(x)C_2(x)}$$
$$\left(\log\left(\frac{\xi_2(x)}{\xi_1(x)}\right)'\right)^2=\frac{-1944 x(E^2\lambda+2450)}{C_1(x)C_2(x)^2}.$$ Moreover, the wronskian of $\xi_1(x)$ and $\xi_2(x)$, which was given in Eq. \eqref{eq:wrons}, is constant.

Thus, we obtain that the differential Galois group is $G(L/K)=\mathbb{D}_\infty$ (the infinite dihedral group).

\subsection*{Acknowledgements}

\noindent PA-H thanks the partial
support  in the final stage of this research to  Ministerio de Educaci\'on Superior, Ciencia y Tecnolog\'ia (MESCYT)  through  grant  2022-1D2-091 of  the  \emph{Fondo  Nacional
de  Innovaci\'on y Desarrollo Cient\'ifico y Tecnol\'ogico (FONDOCYT)}. He also acknowledges to \href{www.spc.edu.do}{Sembrando Pensamiento Cient\'ifico} for the valuable comments of the CEO.\\
JM is member of the Universidad Politécnica de Madrid research group \emph{Modelos Matem\'aticos no lineales} and he thanks to Instituto de Matem\'atica at UASD by its hospitality during his research stay to conclude this paper. \\
T.S. is member of the Physics Institute at UFRJ,


\begin{thebibliography}{99}
	\footnotesize\itemsep=0pt
	
	\bibitem{abst} M. Abramowitz, I. A. Stegun,  Handbook of mathematical functions: with formulas, graphs and mathematical tables, Dover Publications, New York, 1965.
	
	\bibitem{acbl} P. Acosta-Hum\'anez, D. Bl\'azquez-Sanz,  Non-integrability of some hamiltonians with rational potentials,
	\href{http://dx.doi.org/10.3934/dcdsb.2008.10.265}{http://dx.doi.org/10.3934/dcdsb.2008.10.265}, Discrete \& Continuous Dynamical Systems - B, 10 (2008), 265--293, 
	
\bibitem{almp2} P. Acosta-Hum\'anez, J. T.  L\'azaro, J. J. Morales-Ruiz, C. Pantazi,  Differential Galois theory and non-integrability of planar polynomial vector fields, Journal of Differential Equations, 264 (2018), 7183--7212, \href{http://dx.doi.org/10.1016/j.jde.2018.02.016}{http://dx.doi.org/10.1016/j.jde.2018.02.016} 
	
	\bibitem{amw} P. Acosta-Hum\'anez, J. J.  Morales-Ruiz, J.-A. Weil, Galoisian approach to integrability of Schr\"odinger
	equation, Report on Mathematical Physics, 67 (2011), 305--374, \href{http://dx.doi.org/10.1016/S0034-4877(11)60019-0}{http://dx.doi.org/10.1016/S0034-4877(11)60019-0} 

\bibitem{apt} P. B. Acosta-Hum\'anez, M. van der Put, J. Top, Variations for Some Painlev\'e Equation. Symmetry, Integrability and Geometry: Methods and Applications,(2019) 088,10pp. 10pp. \href{ https://doi.org/10.3842/SIGMA.2019.088}{ https://doi.org/10.3842/SIGMA.2019.088}	

\bibitem{ay} P.B. Acosta-Hum\'anez, K. Yagasaki, K. (2020). Nonintegrability of the unfoldings of codimension-two bifurcations,  Nonlinearity, 33 (2020), 1366--1387.

	
	\bibitem{barrow} J. D. Barrow, Chaotic behaviour in general relativity, Physics Reports, 85, No. 1 (1982), 1--49, \href{https://doi.org/10.1016/0370-1573(82)90171-5}{https://doi.org/10.1016/0370-1573(82)90171-5}
	
	
	\bibitem{belinskii2} V. A. Belinskii, E. M. Lifshitz, I. M. Khalatnikov, Oscillatory approach to the singular point in relativistic cosmology, Soviet Physics (Uspekhi), 13 (1972), 745--765.
	
	
	\bibitem{bogo} O. I. Bogoyavlensky, Qualitative Theory of Dynamical Systems in Astrophysics and Gas Dynamics, Springer, Berlin, 1985.
	
	\bibitem{bono} O. I. Bogoyavlensky,  S. P. Novikov, Singularities of the cosmological 
model of the Bianchi IX type according to the qualitative theory of differential equations, Soviet Phys. JETP 37 (1973) 747--755.

\bibitem{crha} T. Crespo, Z. Hajto, Algebraic Groups and Differential Galois Theory, Graduate Studies in Mathematics, 122, American Mathematical Society, Providence, Rhode Island, 2011.



\bibitem{calzeta} E. Calzetta, Homoclinic Chaos in Relativistic Cosmology
 in Deterministic Chaos in General Relativity, 203--235, Springer, New York, 1994, \href{https://doi.org/10.1007/978-1-4757-9993-4_12}{https://doi.org/10.1007/978-1-4757-9993-4-12}
 

\bibitem{cornish}  N. J. Cornish and J. J. Levin,  The mixmaster Universe is chaotic, Phys. Rev. Lett. 78 (1997), 998--1001, \href{https://doi.org/10.1103/PhysRevLett.78.998}{https://doi.org/10.1103/PhysRevLett.78.998}

	

	


\bibitem{matsas} G. Francisco and G. E. A. Matsas, Qualitative and numerical study of Bianchi IX models. General relativity and gravitation, 20 (1998), 1047--1054, \href{https://doi.org/10.1007/BF00759025}{https://doi.org/10.1007/BF00759025}

\bibitem{humpreys} J. Humphreys, Linear algebraic groups, Springer, Berlin, 2012, \href{https://doi.org/10.1007/978-1-4684-9443-3}{https://doi.org/10.1007/978-1-4684-9443-3}
	
	\bibitem{kovacic:1986} J. J. Kovacic, An algorithm for solving second order linear homogeneous differential equations, J. Symbolic Computation 2 (1986), 3--43, \href{https://doi.org/10.1016/S0747-7171(86)80010-4}{https://doi.org/10.1016/S0747-7171(86)80010-4}



\bibitem{misner} C. W. Misner, Mixmaster Universe, Phys. Rev. Lett. 22, (1969), 1071, \href{https://doi.org/10.1103/PhysRevLett.22.1071}{https://doi.org/10.1103/PhysRevLett.22.1071}

\bibitem{mo} J. J. Morales-Ruiz,  Differential Galois Theory and Non-integrability of Hamiltonian Systems, Progress in Mathematics series, 179, Birkh\"ausser, Basel, 1999, \href{https://doi.org/10.1007/978-3-0348-0723-4}{https://doi.org/10.1007/978-3-0348-0723-4}

\bibitem{mora1} J. J. Morales-Ruiz, J.-P. Ramis, Galoisian Obstructions to Integrability of Hamiltonian Systems, Methods and Applications of Analysis, 8, No. 1 (2001), 33--96.

\bibitem{mora2} J. J. Morales-Ruiz, J.-P. Ramis, Galoisian Obstructions to Integrability of Hamiltonian Systems, Methods and Applications of Analysis, 8, No. 1 (2001), 97--112.

\bibitem{morasi} J. J. Morales-Ruiz, J.-P. Ramis, C. Sim\'o, Integrability of Hamiltonian systems and differential Galois groups of higher variational equations, Annales Scientifiques de l'Ecole Normale Sup\'erieure, 40, No. 5 (2007), \href{https://doi.org/10.1016/j.ansens.2007.09.002}{https://doi.org/10.1016/j.ansens.2007.09.002}

\bibitem{oliveira} H. P. de Oliveira, I. D. Soares, T. J. Stuchi, Chaos in anisotropic preinflationary universes, Phys. Rev. D 56 (1997), 730, \href{https://doi.org/10.1103/PhysRevD.56.730}{https://doi.org/10.1103/PhysRevD.56.730}

\bibitem{vasi} M. van der Put, M. Singer, Galois theory of linear differential equations, Grundlehren der mathematischen Wissenschaften, 328,  Springer Verlag, New York 2003, \href{https://doi.org/10.1007/978-3-642-55750-7}{https://doi.org/10.1007/978-3-642-55750-7}
	
\bibitem{ronveaux} A. Ronveaux, Heun's differential equations, Oxford University Press, Oxford, 1995.

\bibitem {stuchi}  I. D. Soares and T. J. Stuchi.  Homoclinic chaos in axisymmetric Bianchi IX cosmologies reexamined, Phys. Rev. D 72, (2005), 083516, \href{https://doi.org/10.1103/PhysRevD.72.083516}{https://doi.org/10.1103/PhysRevD.72.083516}

\bibitem {stuchi1}  I. D. Soares and T. J. Stuchi, Erratum Homoclinic chaos in axisymmetric Bianchi IX cosmologies reexamined, Phys. Rev. D 73, (2006), 069901, \href{https://doi.org/10.1103/PhysRevD.73.069901}{https://doi.org/10.1103/PhysRevD.73.069901}
	
\bibitem{taub} A. H. Taub, Empty Space-Times admitting a Three Parameter Group of motions, Annals of Mathematics, 53 (1951), 472--490, \href{https://doi.org/10.2307/1969567}{https://doi.org/10.2307/1969567}

\end{thebibliography}

\end{document}